\documentclass[11pt,a4paper, final, twoside]{article}
\usepackage{amsmath}
\usepackage{amssymb}
\usepackage{mathrsfs}
\usepackage{mathtools}
\usepackage{fancyhdr}
\usepackage{amsthm}
\usepackage{amsfonts}
\usepackage{amssymb}
\usepackage{amscd}
\usepackage{latexsym}
\usepackage{amsthm}
\usepackage{graphicx}
\usepackage{graphics}
\usepackage{graphtex} 
\usepackage{afterpage}
\usepackage[colorlinks=true, urlcolor=blue,  linkcolor=black, citecolor=black]{hyperref}
\usepackage{color}
\def\bs{\blacksquare}

\newcommand{\HRule}{\rule{\linewidth}{1pt}}
\newcommand{\pf}{\noindent \textit{Proof}:\ }

\newtheorem{thm}{Theorem}[section]

\theoremstyle{proposition}

\theoremstyle{definition}
\newtheorem{defn}{Definition}[section]
\theoremstyle{remark}

\numberwithin{equation}{section}
\begin{document}
	\rhead{\includegraphics[width=14cm]{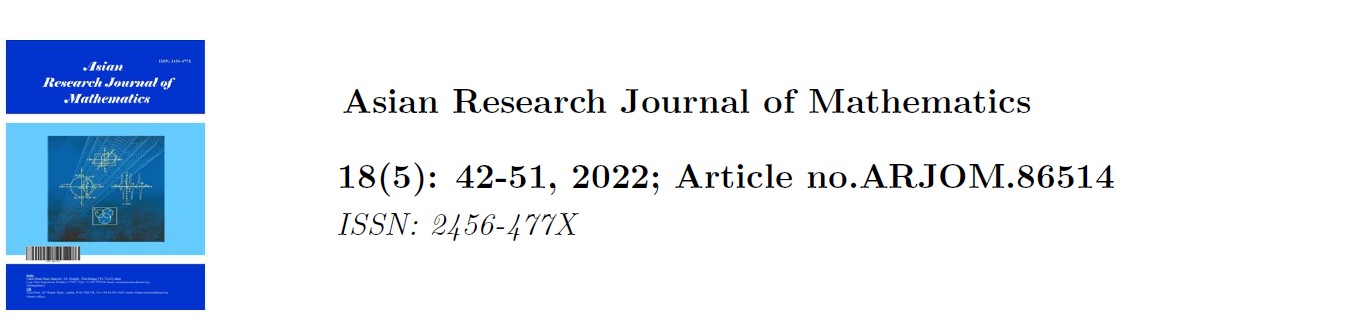}}
\hyphenpenalty=100000
\begin{flushright}
	
{\Large \textbf{\\Harmonic Centrality and Centralization \\of Some Graph Products
}}\\[5mm]
{\large \textbf{Jose Mari E. Ortega$^\mathrm{1^*}$\footnote{\emph{*Corresponding author: E-mail: josemari.ortega@smc.pshs.edu.ph}}, Rolito G. Eballe$^\mathrm{2}$}}\\[3mm]
$^\mathrm{1}${\footnotesize \it Mathematics Unit, Philippine Science High School Southern Mindanao Campus \\ Brgy. Sto. Ni\~no, Davao City, Philippines\\ 
$^\mathrm{2}$Mathematics Department, College of Arts and Sciences, Central Mindanao University\\ Musuan, Maramag, Bukidnon, Philippines}\\ 
\end{flushright}\footnotesize
\indent\quad\indent\quad\indent\quad\indent\quad\indent\quad\indent\quad\indent\quad\indent\quad\indent\indent\quad\quad~~~~~~~~{\it \textbf{Received: 20 February 2022\\ [1mm]
	\indent\quad\indent\quad\indent\quad\indent\quad\indent\quad\indent\quad\indent\quad\indent\quad\indent\quad\quad~~
	{Accepted: 26 April 2022}\\
	\fbox{\large{\bf{\slshape{Original Research Article}}}}~~~~~~~~~~~~~~~~~~~~~~~~~~~\quad Published: 05 May 2022 }}\\[2mm]
\HRule\\[3mm]
{\Large \textbf{Abstract}}\\[4mm]
\fbox{%
\begin{minipage}{5.4in}{\footnotesize Harmonic centrality calculates the importance of a node in a network by adding the inverse of the geodesic distances of this node to all the other nodes. Harmonic centralization, on the other hand, is the graph-level centrality score based on the node-level harmonic centrality. In this paper, we present some results on both the harmonic centrality and harmonic centralization of graphs resulting from some graph products such as Cartesian and direct products of the path $P_2$ with any of the path $P_m$, cycle $C_m$, and fan $F_m$ graphs.

} \end{minipage}}\\[4mm]

\section{Introduction}\label{I1}
Centrality in graph theory and network analysis is based on the importance of a vertex in a graph. Freeman \cite{Fre} tackled the concept of centrality being a salient attribute of social networks which may relate to other properties and processes of said networks. 

Among the many measures of centrality in social network analysis is harmonic centrality. Introduced in 2000 by Marchiori and Latora \cite{Mar}, and developed independently by Dekker \cite{Dek} and Rochat \cite{Roc}, it sums the inverse of the geodesic distances of a particular node to all the other nodes, where it is zero if there is no path from that node to another. It is then normalized by dividing by $m-1$, where $m$ is the number of vertices in the graph. For a related work on harmonic centrality in some graph families, please see \cite{Ort1}. 

While centrality quantifies the importance of a node in a graph, centralization quantifies a graph-level centrality score based on the various centrality scores of the nodes. It sums up the differences in the centrality index of each vertex from the vertex with the highest centrality, with this sum normalized by the most centralized graph of the same order. Centralization may be used to compare how central graphs are. For a related work on harmonic centralization of some graph families, please see \cite{Ort2}. 

In this paper, we derive the harmonic centrality of the vertices and the harmonic centralization of products of some important families of graphs. Graphs considered in this study are the path $P_m$, cycle $C_m$, and fan $F_m$ with the Cartesian and direct product as our binary operations. Our motivation is to be able to obtain formulas that could be of use when one determines the harmonic centrality and harmonic centralization of more complex graphs. For a related study on graph products using betweenness centrality, see \cite{Kum}. 

All graphs under considered here are simple,  undirected, and finite. For graph theoretic terminologies not specifically defined nor described in this study, please refer to \cite{Bal}. 


\section{Preliminary Notes}\label{I2}

For formality, we present below the definitions of the concepts discussed in this paper. 

\begin{defn} \cite{Bol}
	Let $G=(V(G), E(G))$ be a nontrivial graph of order $m$. If $u\in V(G)$,  then the harmonic centrality of $u$ is given by the expression 
	$$ \mathcal{H}_G(u)=\frac{\mathcal{R}_G(u)}{m-1},$$
	where $\mathcal{R}_G(u)=\sum_{x\neq u}\frac{1}{d(u,x)}$ is the sum of the reciprocals of the shortest distance $d(u,x)$ in $G$ between vertices $u$ and $x$, for each $x\in (V(G)\setminus{u})$, with $\frac{1}{d(u,x)}=0$ in case there is no path from $u$ to $x$ in $G$. 
\end{defn} 

\begin{defn} The harmonic centralization of a graph $G$ of order $m$ is given by
	$$ {C_\mathcal{H} (G) } = \frac{\sum\limits_{i=1}^{m}(\mathcal{H}_{G \text{max}}(u)-\mathcal{H}_{G}(u_i))}{\frac{m-2}{m}}   $$
	where $\mathcal{H}_{G \text{max}}(u)$ is the largest harmonic centrality of vertex $u$ in $G$. 
\end{defn} 

In graph $G$ in Figure 1, we have the harmonic centrality of each node calculated as  $\mathcal{H}_{G}(a_1)=\frac{37}{72}, \mathcal{H}_{G}(a_2)=\frac{29}{36}, \mathcal{H}_{G}(a_3)=\frac{23}{36}, \mathcal{H}_{G}(a_4)=\frac{5}{6}, \mathcal{H}_{G}(a_5)=\frac{3}{4}, \mathcal{H}_{G}(a_6)=\frac{13}{18}$, and $\mathcal{H}_{G}(a_7)=\frac{35}{72}$. Clearly, $\mathcal{H}_{G \text{max}}(u)=\frac{5}{6}$ from node $a_4$; thus, the harmonic centralization of graph $G$ has a value of  
${C_\mathcal{H} (G)=\frac{(\frac{5}{6}-\frac{37}{72})+(\frac{5}{6}-\frac{29}{36})+(\frac{5}{6}-\frac{23}{36})+(\frac{5}{6}-\frac{3}{4})+(\frac{5}{6}-\frac{13}{18})+(\frac{5}{6}-\frac{35}{72})}{\frac{7-2}{2}} }= \frac{\frac{13}{12}}{\frac{5}{2}}=\frac{13}{30}$

\begin{figure}[!htb]
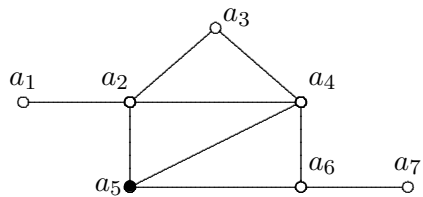

	\Magnify 0.8
	$$\pic 
	\Path (310,60) (310,100) (350,135) (390,100) (390,60) (310,60)
	\Path (310,100) (390,100) 
	\Path (260,100) (310,100) 
	\Path (390,60) (440,60)
	\Path (310,60) (390,100)
	
	\Align [c] ($a_{1}$) (260,110) 
	\Align [c] ($a_{2}$) (303,110) 
	\Align [c] ($a_{3}$) (360,140) 
	\Align [c] ($a_{4}$) (400,110)
	\Align [c] ($a_5$) (300,60)
	\Align [c] ($a_{6}$) (400,70)
	\Align [c] ($a_{7}$) (440,70) 
	\Align [c] ($\bullet$) (310,60)
	\cip$$
	\caption{Graph $G$ with $u\in V(G)$, where $\mathcal{H}_G(a_5)=\frac{3}{4}$ and $C_\mathcal{H}(G)=\frac{13}{30}$}
\end{figure}
\vspace{10pt}

\begin{defn}  
	The $n^{th}$ harmonic number $H_n$ is the sum of the reciprocals of the first $n$ natural numbers, that is
	$ H_n=\sum_{k=1}^n\frac{1}{k}.  $
\end{defn} 

The  definitions for the binary operations of graphs considered in this study are given below. 

\begin{defn} [Cartesian product of graphs]  Let $G$ and $H$ be graphs. The vertex set of the Cartesian product $G \square H$ is the Cartesian product $V(G) \times V(H)$; where two vertices $(u,v)$ and $(u', v')$ are adjacent in $G \square H$ if and only if either (i) $u=u'$ and $v$ is adjacent to $v'$ in $H$, or (ii) $v=v'$ and $u$ is adjacent to $u'$ in $G$. \cite{Bal}
\end{defn} 

\begin{defn} [Direct product of graphs]  Let $G$ and $H$ be graphs. The vertex set of the direct product $G \times H$ is the Cartesian product $V(G) \times V(H)$; where two vertices $(u,v)$ and $(u', v')$ are adjacent in $G \times H$ if and only if (i) $u$ and $u'$ are adjacent in $G$ and (ii) $v$ and $v'$ are adjacent in $H$. \cite{Bal}
\end{defn} 

Drawn below are the paths $P_2$ and $P_3$, and the products $P_2 \square P_3$ and $P_2 \times P_3$. 
\vspace{-5mm}

\begin{center}
	\begin{figure}[!htb]
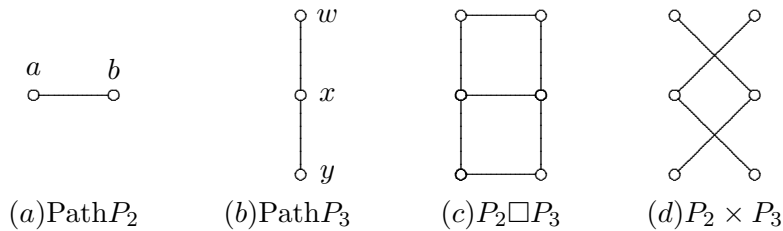

		\Magnify 1.0
		$$\pic 
		\Path (-120,120) (-90,120) 
		\Align [c] ($a$) (-120,130)
		\Align [c] ($b$) (-90,130)
		\Path (-20,150) (-20,120) (-20,90)
		\Align [c] ($w$) (-10,150)
		\Align [c] ($x$) (-10,120)
		\Align [c] ($y$) (-10,90)
		
		\Path (40,90) (40,120) (40,150) (70,150) (70, 120) (70, 90) (40,90) 
		\Path (40,120) (70,120) 
		
		\Path (120,90) (150,120) (120,150) 
		\Path (150,90) (120, 120) (150, 150) 
		
		\Align [c] ($(a)\text{Path} P_2$) (-105,75)
		\Align [c] ($(b)\text{Path} P_3$) (-25,75)
		\Align [c] ($(c)P_2 \square P_3$) (55,75)
		\Align [c] ($(d)P_2\times P_3$) (135,75)
		
		\cip$$
		\caption[The Cartesian Product $P_2 \square P_3$; Direct Product $P_2 \times P_3$.]{(a) Path $P_2$ (b) Path $P_3$ (c) Cartesian Product $P_2 \square P_3$; and (d) Direct Product $P_2 \times P_3$.}
		\label{fig:Products}
	\end{figure}
\end{center}
\vspace{-40pt}

\section{Main Results}\label{J}

The Cartesian and direct products of special graph families considered in this paper are again that of path $P_2$ with path $P_m$ of order $m\geq 1$, path $P_2$ with cycle $C_m$ of order $m\geq 3$, and path $P_2$ with fan $F_m$ of order $m+1\geq 2$.

\begin{center}
	\begin{figure}[!htb]
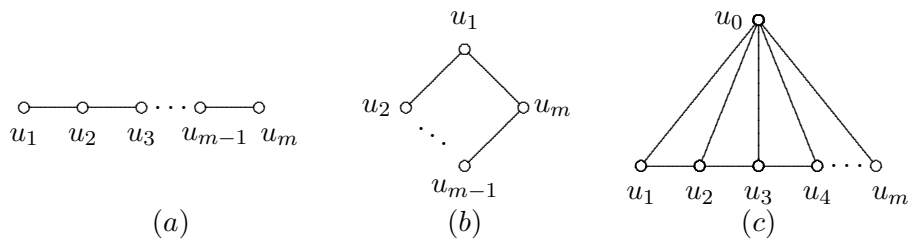

		\Magnify 1.1
		$$\pic 
		\Path (-200,20) (-180,20) (-160,20)
		\Path (-140,20) (-120,20)
		\Align[c] ($\ldots$) (-150,20)
		\Align [c] ($u_{1}$) (-200,10)
		\Align [c] ($u_{2}$) (-180,10)
		\Align [c] ($u_{3}$) (-160,10)
		\Align [c] ($u_{m-1}$) (-135,10)
		\Align [c] ($u_{m}$) (-113,10)
		\Path  (-70,20) (-50,40) (-30,20) (-50,0)
		\Align[c] ($\ddots$) (-61,13)
		\Align [c] ($u_{1}$) (-50,50)
		\Align [c] ($u_{2}$) (-80,20)
		\Align [c] ($u_{m-1}$) (-50,-8)
		\Align [c] ($u_{m}$) (-20,20)
		\Path (50,50) (50,0)
		\Path (50,50) (70,0)
		\Path (50,50) (90,0)
		\Path (50,50) (30,0)
		\Path (50,50) (10,0)
		\Path (10,0) (30,0)
		\Path (50,0) (30,0)
		\Path (70,0) (50,0)
		
		\Align [c] ($u_0$) (40,50)
		\Align [c] ($u_1$) (10,-10)
		\Align [c] ($u_2$) (30,-10)
		\Align [c] ($u_3$) (50,-10)
		\Align [c] ($u_4$) (70,-10)
		\Align [c] ($u_m$) (95,-10)
		\Align [c] ($\ldots$) (80,0)
		
		\Align [c] ($(a)$) (-150,-20)
		\Align [c] ($(b)$) (-50,-20)
		\Align [c] ($(c)$) (50,-20)
		\cip$$
		\caption[The Path $P_m$;\; Cycle $C_m$;\; and Fan $F_m$.]{(a) Path $P_m$;\;(b) Cycle $C_m$;\; and (c) Fan $F_m$.}
		\label{fig:Path}
	\end{figure}
\end{center}

\begin{thm} 
	\normalfont
	For the Cartesian product of the path $P_2=[u_1, u_2]$ and the path $P_m=[v_1, v_2, ..., v_m]$, the harmonic centrality of any vertex $(u_i, v_j)$ is given by
	
	\[
	\mathcal{H}_{P_2\square P_m}(u_i, v_j) =  \begin{dcases*}
		\frac{1}{2m-1} \Big(2H_{m-1}+\frac{1}{m}\Big) &  \text{for $1\leq i \leq 2$, $j=1$ or $j=m$} \\
		\frac{1}{2m-1} \Big[2\Big(H_{j-1}+H_{m-j}\Big) &  \\
		\quad \quad + \frac{1}{j} +\frac{1}{m-j+1}-1\Big] & \text{for $1\leq i \leq 2$, $1<j<m$.}  \\
	\end{dcases*}
	\]
	
	\pf The product of $P_2$ and $P_m$ is also called a ladder graph $L_m$ of order $2m$. Considering the structure of the product $P_2\square P_m$ or $L_m$, we can partition its vertex set into two subsets  $V_2(L_m)=\{(u_1, v_1), (u_1, v_2),..., (u_1, v_m)\}$ and $V_2(L_m)=\{(u_2, v_1), (u_2, v_2),..., (u_2, v_m)\}$ with $P_2=[u_1, u_2]$ and $P_m=[v_1, v_2, ..., v_m]$. For $(u_i, v_1)$ and $(u_i, v_m)$ with $i=1, 2$ we have \\
	
	\noindent \small 
	$\begin{aligned}
		\mathcal{R}_{L_m}(u_i,v_1)	&  = \mathcal{R}_{L_m}(u_i,v_m) \ = \mathcal{R}_{L_m}(u_1,v_1) \ \\
		& = \sum\limits_{x \in V_1(L_{m}),x\neq (u_1, v_1)} \frac{1}{\text{d}_{L_{m}} ((u_1,v_1), x)} +\sum\limits_{x \in V_2(L_{m})} \frac{1}{\text{d}_{L_{m}} ((u_1,v_1), x)}  \\
		& = \frac{1}{\text{d}((u_1,v_1), (u_1,v_2))}+\frac{1}{\text{d} ((u_1,v_1), (u_1,v_3))}+\\
		& \quad \quad ...+ \frac{1}{\text{d}((u_1,v_1), (u_1,v_m))}+\frac{1}{\text{d} ((u_1,v_1), (u_2,v_1))}+\\
		& \quad \quad ...+\frac{1}{\text{d} ((u_1,v_1), (u_2,v_{m-1}))}+\frac{1}{\text{d} ((u_1,v_1), (u_2,v_m)}\\
		& = \Big[1+\frac{1}{2}+...+\frac{1}{m-1}\Big]+ \Big[1+\frac{1}{2}+...+\frac{1}{m-1}+\frac{1}{m}\Big] \\
		& = 2\sum_{k=1}^{m-1} \frac{1}{k} + \frac{1}{m} \\
		& = 2H_{m-1}+\frac{1}{m} 
	\end{aligned}$ \\
	
	\noindent As for $(u_i, v_j)\in V(L_m)$, where $i=1, 2,$ and $j=2, 3, ..., m-1$, 
	
	\noindent \small
	$\begin{aligned}
		\mathcal{R}_{L_m}(u_i, v_j)\ & =\mathcal{R}_{L_m}(u_1, v_j)\  \\
		& = \sum\limits_{x \in V_1(L_{m}),x\neq (u_1, v_j)} \frac{1}{\text{d}_{L_{m}} ((u_1,v_j), x)} +\sum\limits_{x \in V_2(L_{m})} \frac{1}{\text{d}_{L_{m}} ((u_1,v_j), x)}  \\
		& = \sum\limits_{1\leq k\leq m,k\neq j} \frac{1}{\text{d}_{L_{m}} ((u_1,v_j), (u_1, v_k))} +\sum\limits_{k=1}^{m} \frac{1}{\text{d}_{L_{m}} ((u_1,v_j), (u_2, v_k))}  \\
		& = \big[\frac{1}{j-1}+\frac{1}{j-2}+...+\frac{1}{3} +\frac{1}{2}+1\big]+\big[1+\frac{1}{2}+\frac{1}{3}+...+\frac{1}{m-j}\big]\\
		& \quad \quad +\big[\frac{1}{j}+\frac{1}{j-1}+...+\frac{1}{3} +\frac{1}{2}\big]+\big[1+\frac{1}{2}+\frac{1}{3}+...+\frac{1}{m-j+1}\big]\\
		& = \big[\frac{1}{j-1}+\frac{1}{j-2}+...+\frac{1}{3} 	+\frac{1}{2}+1\big]+\big[1+\frac{1}{2}+\frac{1}{3}+...+\frac{1}{m-j}\big]\\
		& \quad \quad +\big[\frac{1}{j-1}+\frac{1}{j-2}+...+\frac{1}{3} +\frac{1}{2}+1\big]+\big[1+\frac{1}{2}+\frac{1}{3}+...+\frac{1}{m-j}\big] \\
		& \quad \quad \quad + \frac{1}{j}+\frac{1}{m-j+1}-1\\
		& = 2(H_{j-1}+H_{m-j}) + \frac{1}{j} +\frac{1}{m-j+1}-1
	\end{aligned}$

	Thus, after normalizing we get
	\[
	\mathcal{H}_{P_2\square P_m}(u_i, v_j) =  \begin{dcases*}
		\frac{1}{2m-1} \Big(2H_{m-1}+\frac{1}{m}\Big) &  \text{for $1\leq i \leq 2$, $j=1$ or $j=m$} \\
		\frac{1}{2m-1} \Big[2\Big(H_{j-1}+H_{m-j}\Big) &  \\
		\quad \quad + \frac{1}{j} +\frac{1}{m-j+1}-1\Big] & \text{for $1\leq i \leq 2$, $1<j<m$.}  \\
	\end{dcases*}
	\]\\\\ [-40pt]
\end{thm}

\begin{thm} 
	\normalfont 
	For the Cartesian product of path $P_2=[u_1, u_2]$ of order 2 and cycle graph $C_m=[v_1, v_2, ..., v_m, v_1]$ of order $m$, the harmonic centrality of any vertex $(u_i, v_j)$ is given by 
	\[
	\mathcal{H}_{P_2 \square C_m}(u_i,v_j) =  \begin{dcases*}
		\frac{1}{2m-1} \Big(4H_\frac{m-1}{2}+\frac{3-m}{m+1}\Big) & \text{if $m$ is odd}\\
		\frac{1}{2m-1} \Big(4H_{\frac{m}{2}}+\frac{2}{m+2}-\frac{m+2}{m}\Big)  & \text{if $m$ is even}.
	\end{dcases*}
	\] 
	
	\begin{sloppypar}		
		\pf The Cartesian product of $P_2$ and $C_m$ is also known as a prism $Y_m$ of order $2m$. Considering its structure, we can partition its into two subsets $V_1(Y_m)={(u_1, v_1), (u_1, v_2),..., (u_1, v_m)}$ and $V_2(Y_m)={(u_2, v_1), (u_2, v_2),..., (u_2, v_m)}$. If $m$ is odd, we have 
	\end{sloppypar}
	
	\noindent
	$\begin{aligned} 
		\mathcal{R}_{Y_m}(u_i,v_j)  \ & =  \mathcal{R}_{Y_m}(u_1,v_j)  \	\\
		& = \sum\limits_{x \in V_1(Y_{m}), x\neq (u_1, v_j)} \frac{1}{\text{d}_{Y_m} ((u_i,v_j), x)} +\sum\limits_{x \in V_2(Y_m)} \frac{1}{\text{d}_{Y_m} ((u_i,v_j), x)}  \\
		& = \Big[1+1+...+\frac{1}{\frac{m-1}{2}}+\frac{1}{\frac{m-1}{2}}\Big]+ \Big[1+\frac{1}{2}+\frac{1}{2}+...+\frac{1}{\frac{m+1}{2}}+\frac{1}{\frac{m+1}{2}}\Big] \\
		& = 4\sum_{k=1}^{\frac{m-1}{2}} \frac{1}{k} + 2\Big(\frac{2}{m+1}\Big) -1 \\
		& = 4H_{\frac{m-1}{2}}+\frac{3-m}{m+1} 
	\end{aligned}$ \\
	
	\noindent
	If $m$ is even, we have
	
	\noindent
	$\begin{aligned} 
		\mathcal{R}_{Y_m}(u_i,v_j)  \	& = \sum\limits_{x \in V_1(Y_{m})} \frac{1}{\text{d}_{Y_m} ((u_i,v_j), x)} +\sum\limits_{x \in V_2(Y_m)} \frac{1}{\text{d}_{Y_m} ((u_i,v_j), x)}  \\
		& = \Big[1+1+\frac{1}{2}+\frac{1}{2}+...+\frac{1}{\frac{m}{2}}\Big]+ \Big[1+\frac{1}{2}+\frac{1}{2}+...+\frac{1}{\frac{m+2}{2}}\Big]\\
		& = 4\sum_{k=1}^{\frac{m}{2}} \frac{1}{k}  + \frac{1}{\frac{m+2}{2}} -1 -\frac{1}{\frac{m}{2}}  \\
		& = 4H_{\frac{m}{2}}+\frac{2}{m+2}-\frac{m+2}{m}\\ 
	\end{aligned}$ \\
	
	Normalizing and consolidating these results we get 
	\[
	\mathcal{H}_{Y_m}(u_i,v_j) =  \begin{dcases*}
		\frac{1}{2m-1} \Big(4H_\frac{m-1}{2}+\frac{3-m}{m+1}\Big) & \text{if $m$ is odd} \\
		\frac{1}{2m-1} \Big(4H_{\frac{m}{2}}+\frac{2}{m+2}-\frac{m+2}{m}\Big)  & \text{if $m$ is even.} $\quad \bs$ 
	\end{dcases*}
	\] \\
	
\end{thm}

\begin{thm} 
	\normalfont 
	For the Cartesian product of path $P_2$ of order 2 and fan graph $F_m$ of order $m+1$ in Figure 3c, the harmonic centrality of any vertex $(u_i, v_j)$ is given by
	
	\[
	\mathcal{H}_{P_2\square F_m }(u_i, v_j) =  \begin{dcases*}
		\frac{3m+2}{2(2m+1)}  & \text{for $i=1,2$ and $j=0$} \\
		\frac{5m+14}{6(2m+1)}  &  \text{for $i=1,2$ and $j=1,m$} \\
		\frac{5m+18}{6(2m+1)} & \text{for $i=1,2$ and $1<j<m$} \\
	\end{dcases*}
	\]
	\pf Vertices $(u_1, j_0)$ and $(u_2, j_0)$ are adjacent to $m+1$ vertices, and have a distance of 2 to $m$ vertices. Thus, 
	
	$$\mathcal{H}_{P_2\square F_m}(u_i, v_0)\ = \frac{1}{2m+1}\Big[1(m+1) + \frac{1}{2} (m)\Big] = \frac{3m+2}{2(2m+1)} $$ 
	
	For vertices $(u_1, v_1), (u_1, v_m), (u_2, v_1)$ and $(u_2, v_m)$, each is adjacent to three vertices, each has a distance of 2 to $m$ vertices, and each has a distance of 3 to $m-2$ vertices; therefore,
	
	$$\mathcal{H}_{P_2\square F_m}(u_i, v_1)\ = \mathcal{H}_{P_2\square F_m}(u_i, v_m)\ = \frac{1}{2m+1}\Big[1(3) + \frac{1}{2} (m) + \frac{1}{3} (m-2)\Big]  = \frac{5m+14}{6(2m+1)} $$  
	
	As for vertices $(u_i, v_j)$ for $i=1,2$ and $1<j<m$, each is adjacent to four vertices, each has a distance of 2 to $m$ vertices, and each has a distance of 3 to $m-3$ vertices; thus, 
	
	$$\mathcal{H}_{P_2\square F_m}(u_i, v_j)\	 = \frac{1}{2m+1}\Big[1(4) + \frac{1}{2} (m) + \frac{1}{3} (m-3)\Big]  = \frac{5m+18}{6(2m+1)}. \quad \bs  $$  	
\end{thm}

\begin{thm} 
	\normalfont For the Cartesian product of $P_2=[u_1, u_2]$ of order 2 and path $P_m=[u_1, u_2, ..., u_m]$ of order $m$, the harmonic centralization is given by
	
	$$ {\footnotesize
		{C}_\mathcal{H}(L_m) = \begin{dcases*}
			\frac{4}{(m-1)(2m-1)}\Big[2(m-1)H_\frac{m-1}{2}-2H_{m-1} \\ 
			\quad +\frac{2(m-1)}{m+1}-\frac{m-1}{2}-\frac{1}{m} & \text{if $m$ is odd  }  \\ 
			\quad \quad -\sum\limits_{i=2}^{\frac{m-1}{2}}\Big(2H_{i-1}+2H_{m-i}+\frac{1-i}{i}+\frac{1}{m-i+1}\Big)\Big] \\ \\[-8pt]
			\frac{2}{(2m-1)(m-1)}\Big[4(m-2)H_\frac{m}{2}-4H_{m-1}  \\
			\quad -\frac{m^2-2}{m}+\frac{2m-4}{m+2}  & \text{if $m$ is even.} \\ 
			\quad -2\sum\limits_{i=2}^{\frac{m-2}{2}} \Big(2H_{i-1}+2H_{m-i}+ \frac{1-i}{i}+ \frac{1}{m-i+1}\Big)\Big]
		\end{dcases*} 
	}$$ 
	
	\pf The product of $P_2=[u_1, u_2]$  and $P_m=[u_1, u_2, ..., u_m]$ of order is a ladder graph $L_m$ of order $2m$. In a ladder graph, if $m$ is odd, then vertices $(u_1, v_{\frac{m+1}{2}})$ and $(u_2, v_{\frac{m+1}{2}})$ will have the maximum harmonic centrality of $\frac{4}{2m-1}\Big(H_\frac{m-1}{2}+\frac{1}{m+1}-\frac{1}{4}\Big)$. So, 
	
	$$ {\small
		\begin{aligned} 
			{C}_\mathcal{H}(L_m) & = \frac{1}{\frac{2m-2}{2}}\Big[\Big(\frac{4(2(m-1))}{2m-1}\Big(H_\frac{m-1}{2}+\frac{1}{m+1}-\frac{1}{4}\Big)\Big)-\frac{4}{2m-1}\Big(2H_{m-1}+\frac{1}{m}\Big) \\
			& \quad \quad -\frac{4}{2m-1}\sum\limits_{j=2}^{\frac{m-1}{2}}\Big(2H_{j-1}+2H_{m-j}+\frac{1-j}{j}+\frac{1}{m-j+1}\Big)\Big]\\
			& = \frac{4}{(m-1)(2m-1)}\Big[2(m-1)H_\frac{m-1}{2}-2H_{m-1}+\frac{2(m-1)}{m+1}-\frac{m-1}{2}-\frac{1}{m}  \\  	
			& \quad \quad -\sum\limits_{j=2}^{\frac{m-1}{2}}\Big(2H_{j-1}+2H_{m-j}+\frac{1-j}{j}+\frac{1}{m-j+1}\Big)\Big].\\
		\end{aligned}
	}$$ 
	
	On the other hand, if $m$ is even in the product, vertices $(u_1, v_{\frac{m}{2}}), (u_2, v_{\frac{m}{2}}), (u_1, v_{\frac{m+2}{2}})$ and $(u_2, v_\frac{m+2}{2})$ will have the maximum harmonic centrality of $\frac{1}{2m-1}\Big(4H_\frac{m}{2}-\frac{m+2}{m}+\frac{2}{m+2}\Big)$. So, 
	
	$$ {\small
		\begin{aligned} 
			{C}_\mathcal{H}(L_m) & = \frac{1}{\frac{2m-2}{2}}\Big[\frac{4}{2m-1}\Big(4H_\frac{m-1}{2}-2H_{m-1}-\frac{m+3}{m}+\frac{2}{m+2}\Big)  \\  	
			& \quad \quad +\frac{2(m-4)}{2m-1}\Big(4H_\frac{m}{2}-\frac{m+2}{m}+\frac{2}{m+2}\Big) \\
			& \quad \quad -\frac{4}{2m-1}\sum\limits_{j=2}^{\frac{m-2}{2}}\Big(2H_{j-1}+2H_{m-j}+\frac{1}{j}+\frac{1}{m-j+1}-1\Big)\Big]\\
			& = \frac{2}{(m-1)(2m-1)}\Big[8H_\frac{m}{2}+4(m-4)H_\frac{m}{2}-4H_{m-1} \\  	
			& \quad \quad -\frac{2(m+3)-(m+2)(m-4)}{m}+\frac{4+2(m-4)}{m+2} \\
			& \quad \quad -2\sum\limits_{j=2}^{\frac{m-2}{2}}\Big(2H_{j-1}+2H_{m-j}+\frac{1-j}{j}+\frac{1}{m-j+1}\Big)\Big]\\
			& =\frac{1}{(2m-1)(m-1)}\Big[4(m-2)H_\frac{m}{2}-4H_{m-1}  \\ 
			& \quad -\frac{(m^2-2)}{m}+\frac{2m-4}{m+2} \\ 
			& \quad -2\sum\limits_{j=2}^{\frac{m-2}{2}} \Big(2H_{j-1}+2H_{m-j}+ \frac{1-j}{j}+ \frac{1}{m-j+1}\Big)\Big]. \quad \bs  
		\end{aligned}
	}$$
\end{thm} 

\begin{thm} 
	\normalfont
	For the Cartesian product of path $P_2$ of order 2 and Cycle $C_m$ of order $m$, the harmonic centralization is zero. \\
	
	\pf Each vertex of the resulting graph of $P_2\square C_m$  will have the same harmonic centralities of $\mathcal{H}_{P_2\square C_m}(u_i, v_j)=\frac{2}{m-1}(4H_\frac{m-1}{2}+\frac{3-m}{m+1})$ if $m$ is odd and $\mathcal{H}_{P_2\square C_m}(u_i, v_j)= \frac{2}{m-1}(4H_\frac{m-1}{2}+\frac{2}{m+2}-\frac{m+2}{m})$ if $m$ is even. Therefore, this results to a harmonic centralization of zero. $\bs$	
	
\end{thm}

\begin{thm} 
	\normalfont 
	For the Cartesian product of path $P_2$ of order 2 and fan graph $F_m$ of order $m+1$, the harmonic centralization is given by
	\[
	{C}_\mathcal{H}(P_2\square F_m) = \frac{4(m-1)(m-2)}{3m(2m+1)}
	\]
	\pf The vertex with the maximum harmonic centrality would be vertices $(u_1, v_0)$ and $(u_2, v_0)$ having a value of $\frac{3m+2}{2(2m+1)}$. The harmonic centrality of all the other vertices will be subtracted from this maximum harmonic centrality, thus
	
	$\begin{aligned} 
		{C}_\mathcal{H}(P_2\square F_m)\	& = \Big(\frac{1}{\frac{2(m+1)-2}{2}}\Big)\Big[4\Big(\frac{3m+2}{4m+2}-\frac{5m+14}{6(2m+1)}\Big) \\
		& \quad \quad \quad \quad +(2m-4)\Big(\frac{3m+2}{4m+2}-\frac{5m+18}{6(2m+1)}\Big)\Big] \\  	
		& =  \cfrac{4(m-1)(m-2)}{3m(2m+1)}. \quad \bs \\	
	\end{aligned}$ \\ 
	
\end{thm}

\begin{thm} 
	\normalfont 
	For the direct product of path $P_2$ of order 2 and a path graph $P_m$ of order $m$, the harmonic centrality of any vertex $(u_i, v_j)$ is given by
	
	\[
	\mathcal{H}_{P_2\times P_m }(u_i, v_j) =  \begin{dcases*}
		\frac{H_{m-1}}{2m-1}  & \text{for $i=1,2$ and $j=1$ or $m$} \\
		\frac{H_{j-1}+H_{m-j}}{2m-1}  &  \text{for $i=1,2$ and $1<j<m$} \\
	\end{dcases*}
	\]
	\pf The resulting direct product of $P_2\times P_m$ is a graph of order $2m$ composed of a pair of disconnected path graphs of order $m$. Thus, the endpoints will have a harmonic centrality of $\frac{H_{m-1}}{2m-1}$, while the inside vertices will have a harmonic centrality value of $\frac{H_{j-1}+H_{m-j}}{2m-1}$. $\quad \bs$
	
\end{thm}

\begin{thm} 
	\normalfont 
	For the direct product of path $P_2$ of order 2 and a cycle graph $C_m$ of order $m$, the harmonic centrality of any vertex $(u_i, v_j)$ is given by
	\[
	\mathcal{H}_{P_2\times C_m }(u_i, v_j) =  \begin{dcases*}
		\frac{1}{2m-1}\big(2H_{m-1}+\frac{1}{m}\big)  & \text{if $m$ is odd} \\
		\frac{1}{2m-1}\big(2H_{\frac{m-2}{2}}+\frac{2}{m}\big)  &  \text{if $m$ is even} \\
	\end{dcases*}
	\]
	\pf The resulting graph of $P_2\times C_m$ of order $2m$ is a pair of cycle graphs of order $m$. From what we know of cycle graphs $\mathcal{R}_{C_m}(u)=2H_{m-1}+\frac{1}{m}$ if $m$ is odd, while $\mathcal{R}_{C_m}(u)=2H_{\frac{m-2}{2}}+\frac{2}{m}$ if $m$ is even. $\quad \bs$
	
\end{thm}

\begin{thm} 
	\normalfont 
	For the direct product of path $P_2$ of order 2 and a fan graph $F_m$ of order $m+1$, the harmonic centrality of any vertex $(u_i, v_j)$ is given by
	\[
	\mathcal{H}_{P_2\times F_m }(u_i, v_j) =  \begin{dcases*}
		\frac{m}{2m+1}  & \text{for $i=1,2$ and $j=0$} \\
		\frac{m+2}{2(2m+1)}  &  \text{for $i=1,2$ and $j=1$ or $m$} \\
		\frac{m+3}{2(2m+1)}  &  \text{for $i=1,2$ and $1<j<m$} \\
	\end{dcases*}
	\]
	\pf The resulting graph of $P_2\times F_m$ is of order $2m+2$ composed of a pair of fan graphs each with an order of $m+1$. So vertices $(u_1, v_0)$ and $(u_2, v_0)$ will be adjacent to $m$ vertices and normalized by $2m+1$. While vertices $(u_1, v_1)$ and $(u_1, v_m)$ will be adjacent to 2 vertices and have a distance of 2 to $m-2$ vertices. Therefore, $\frac{1}{2m+1}\Big(\frac{m-2}{2}+2\Big) = \frac{m+2}{2(2m+1)}$. All the other vertices will be adjacent to 3 other vertices and have a distance of 2 to $m-3$ vertices, therefore, $\frac{1}{2m+1}\Big(\frac{m-3}{2}+3\Big) = \frac{m+3}{2(2m+1)}$.
	$\quad \bs$
	
\end{thm}

\begin{thm} 
	\normalfont 
	For the direct product of path $P_2$ of order 2 and a path graph $P_m$ of order $m$, the harmonic centralization is given by
	
	\[
	\footnotesize 
	{C}_\mathcal{H}(P_2 \times P_m) =  \begin{dcases*}
		\frac{4\big((m-1)H_{\frac{m-1}{2}}-H_{m-1}-\sum\limits_{j=2}^{\frac{m-1}{2}}(H_{j-1}+H_{m-j})\big)}{(m-1)(2m-1)}  & \text{if $m$ is odd} \\
		\frac{4\big((m-2)(H_{\frac{m-2}{2}}+\frac{1}{m})-H_{m-1}-\sum\limits_{j=2}^{\frac{m-2}{2}}(H_{j-1}+H_{m-j})\big)}{(m-1)(2m-1)}  &  \text{if $m$ is even} \\
	\end{dcases*}
	\]
	\pf If $m$ is odd, then the maximum harmonic centrality will be $\frac{H_{\frac{m-1}{2}}}{2m-1}$ for vertices  $(u_1, v_{\frac{m-1}{2}})$ and $(u_2, v_{\frac{m-1}{2}})$. Other vertices in the graph are the endpoints which has a harmonic centrality of $\frac{H_{m-1}}{2m-1}$ and all other vertices have harmonic centrality of $\frac{H_{j-1}+H_{m-j}}{2m-1}$. Subtracting from the maximum and normalizing we get 	$$\frac{4\big((m-1)H_{\frac{m-1}{2}}-H_{m-1}-\sum\limits_{j=2}^{\frac{m-1}{2}}(H_{j-1}+H_{m-j})\big)}{(m-1)(2m-1)}$$. \\
	On the other hand, if $m$ is even, then vertices $(u_1, v_{\frac{m}{2}}), (u_2, v_{\frac{m}{2}}), (u_1, v_{\frac{m+1}{2}})$ and $(u_2, v_{\frac{m+1}{2}})$ will have the maximum harmonic centrality of $\frac{H_{\frac{m-2}{2}}+\frac{1}{m}}{2m-1}$. The harmonic centrality of all other vertices will be subtracted from this value and normalized to arrive at the harmonic centralization value of  $$\frac{4\big((m-2)(H_{\frac{m-2}{2}}+\frac{1}{m})-H_{m-1}-\sum\limits_{j=2}^{\frac{m-2}{2}}(H_{j-1}+H_{m-j})\big)}{(m-1)(2m-1)} \quad \bs$$.

\end{thm}

\begin{thm} 
	\normalfont 
	For the direct product of path $P_2$ of order 2 and a cycle graph $C_m$ of order $m$, the harmonic centralization is zero. \\
	\pf Since the harmonic centrality of the vertices in the direct product of path $P_2$ and cycle graph $C_m$ have the same values, the harmonic centralization therefore equates to zero.  
	$\quad \bs$
\end{thm}

\begin{thm} 
	\normalfont 
	For the direct product of Path $P_2$ of order 2 and Fan graph $F_m$ of order $m+1$, the harmonic centralization is given by
	\[
	{C}_\mathcal{H}(P_2\times F_m) = \frac{(m-1)(m-2)}{m(2m+1)}
	\]
	\pf The maximum harmonic centrality value is $\frac{m}{2m+1}$ while four vertices have harmonic centrality values of $\frac{m+2}{2(2m+1)}$. As for the other vertices, they have a harmonic centrality of $\frac{m+3}{2(2m+1)}$, so dividing by $m$ to normalize, we have \\
	$\begin{aligned} 
		{C}_\mathcal{H}(P_2\times F_m)\	& = \frac{1}{m}\Big[2m\Big(\frac{m}{2m+1}\Big)-4\Big(\frac{m+2}{2(2m+1)}\Big)-(2m-4)\Big(\frac{m+3}{2(2m+1)}\Big)\Big] \\  	
		& =  \cfrac{(m-2)(m-1)}{m(2m+1)}. \quad \bs \\	
	\end{aligned}$ \\ 
\end{thm}

\section{Conclusions}
Harmonic centrality is one of the more recent centrality measures that identifies the importance of a node, while harmonic centralization quantifies how centralized a graph is based on the node-level harmonic centrality. In this paper, we introduced some results on the harmonic centrality of the nodes and harmonic centralization of graphs resulting from the Cartesian and direct products of the path $P_2$ with any of the path $P_m$, cycle $C_m$, and fan $F_m$ graphs. For further studies, results can be derived for other families of graphs and other binary operations. \\

\noindent \Large\textbf{Acknowledgement}\\[2mm] 
\footnotesize The authors would like to acknowledge the valuable comments and inputs made by the anonymous referees.\\[2mm]

\noindent{\Large\bf Competing Interests}\\\\
Authors have declared that no competing interests exist.

\scriptsize\-----------------------------------------------------------------------------------------------------------------------------------------\\\copyright \it 20YY  Author name; This is an Open Access article distributed under the terms of the Creative Commons Attribution License
\href{http://creativecommons.org/licenses/by/2.0}{http://creativecommons.org/licenses/by/4.0},  which permits unrestricted use, distribution, and reproduction in any medium,
provided the original work is properly cited.
\end{document}